\documentclass{article}[12pt]

\usepackage{latexsym}
\usepackage{amsmath}
\usepackage{amssymb}

\begin{document}

\title{Solution of the system of two coupled first-order ODEs with
second-degree polynomial right-hand sides}

\author{Francesco Calogero$^{a,b}$\thanks{e-mail: francesco.calogero@roma1.infn.it}
\thanks{e-mail: francesco.calogero@uniroma1.it}
 , Farrin Payandeh$^c$\thanks{e-mail: farrinpayandeh@yahoo.com}
 \thanks{e-mail: f$\_$payandeh@pnu.ac.ir}}

\maketitle   \centerline{\it $^{a}$Physics Department, University of
Rome "La Sapienza", Rome, Italy}

\maketitle   \centerline{\it $^{b}$INFN, Sezione di Roma 1}

\maketitle

\maketitle   \centerline{\it $^{c}$Department of Physics, Payame
Noor University, PO BOX 19395-3697 Tehran, Iran}

\maketitle

\begin{abstract}

The \textit{explicit} solution $x_{n}\left( t\right) ,$ $n=1,2,$ of the
\textit{initial-values} problem is exhibited of a \textit{subclass} of the
\textit{autonomous} system of $2$ coupled \textit{first-order} ODEs with
\textit{second-degree} polynomial right-hand sides, hence featuring $12$
\textit{a priori arbitrary} (time-independent) coefficients:%
\begin{equation*}
\dot{x}_{n}=c_{n1}\left( x_{1}\right) ^{2}+c_{n2}x_{1}x_{2}+c_{n3}\left(
x_{2}\right) ^{2}+c_{n4}x_{1}+c_{n5}x_{2}+c_{n6}~,~~~n=1,2~.
\end{equation*}%
The solution is \textit{explicitly} provided if the $12$ coefficients $%
c_{nj} $ ($n=1,2;~j=1,2,3,4,5,6$) are expressed by \textit{explicitly}
provided formulas in terms of $10$ \textit{a priori arbitrary} parameters;
the \textit{inverse} problem to express these $10$ parameters in terms of
the $12 $ coefficients $c_{nj}$ is also \textit{explicitly} solved, but it
is found to imply---as it were, \textit{a posteriori}---that the $12$
coefficients $c_{nj}$ must then satisfy $4$ \textit{algebraic constraints},
which are \textit{explicitly} exhibited. Special subcases are also
identified the \textit{general} solutions of which are \textit{completely
periodic} with a period independent of the initial data ("isochrony"), or
are characterized by additional restrictions on the coefficients $c_{nj}$
which identify particularly interesting models.
\end{abstract}

\section{Introduction}

The system characterized by the $2$ \textit{nonlinearly-coupled} ODEs
\begin{subequations}
\label{1}
\begin{eqnarray}
\dot{x}_{n}\left( t\right) =c_{n1}\left[ x_{1}\left( t\right) \right]
^{2}+c_{n2}x_{1}\left( t\right) x_{2}\left( t\right) +c_{n3}\left[
x_{2}\left( t\right) \right] ^{2} &&  \notag \\
+c_{n4}x_{1}\left( t\right) +c_{n5}x_{2}\left( t\right) +c_{n6}~,~~~n=1,2~,
&&  \label{1a}
\end{eqnarray}%
namely%
\begin{equation}
\dot{x}_{1}\left( t\right) =c_{11}\left[ x_{1}\left( t\right) \right]
^{2}+c_{12}x_{1}\left( t\right) x_{2}\left( t\right) +c_{13}\left[
x_{2}\left( t\right) \right] ^{2}+c_{14}x_{1}\left( t\right)
+c_{15}x_{2}\left( t\right) +c_{16}~,  \label{1b}
\end{equation}%
\begin{equation}
\dot{x}_{2}\left( t\right) =c_{21}\left[ x_{1}\left( t\right) \right]
^{2}+c_{22}x_{1}\left( t\right) x_{2}\left( t\right) +c_{23}\left[
x_{2}\left( t\right) \right] ^{2}+c_{24}x_{1}\left( t\right)
+c_{25}x_{2}\left( t\right) +c_{26}~,  \label{1c}
\end{equation}%
is a prototypical example of \textit{autonomous} dynamical systems. It
features the $12$ \textit{a priori arbitrary} time-independent parameters $%
c_{nj}$ ($n=1,2;$ $j=1,2,3,4,5,6$). The main finding of the present paper is
to show that the \textit{initial-values} problem of this dynamical system
can be \textit{explicitly} solved, provided the $12$ coefficients $c_{nj}$ ($%
n=1,2$; $j=1,2,3,4,5,6$) satisfy $4$ \textit{simple constraints}:\textit{\ }%
for a neat version of these \textit{constraints} see below \textbf{Section 4
}where the solution of the \textit{initial-values} problem of the system (%
\ref{1}) is displayed.

\textbf{Notation 1-1}. The $2$ (possibly \textit{complex}) numbers $%
x_{n}\equiv x_{n}\left( t\right) ,$ $n=1,2,$ are the dependent variables; $t$
is the independent variable ("time"; but the treatment remains valid when $t$
is considered a \textit{complex} number); superimposed dots denote $t$%
-differentiations; the $12$ time-independent (possibly \textit{complex})
numbers $c_{nj},$ $n=1,2,$ $j=1,2,3,4,5,6,$ are parameters. Hereafter the
time-dependence of variables is often not \textit{explicitly} indicated,
when this omission is unlikely to cause misunderstandings. The indices $n,$ $%
m,$ $j,$ $k,$ $\ell $---as indeed generally clear from the context---run
respectively over the integers from $1$ to $2$ ($n=1,2$), from $1$ to $2$ ($%
m=1,2$), from $1$ to $6$ ($j=1,2,3,4,5,6$), from $1$ to $3$ ($k=1,2,3$) and
from $2$ to $0$ ($\ell =2,1,0$). $\blacksquare $

\textbf{Remark 1-1}. The system (\ref{1}) has been investigated over time in
an enormous number of mainly \textit{mathematical}, or mainly \textit{%
applicable}, contexts; too many to allow any attempt to provide a list of
references that would do justice to the multitude of relevant papers. We
limit ourselves here to quote only $3$ very recent papers: \cite{CCL2020}
(where the case with \textit{homogeneous} right-hand sides has been treated,
i. e. the system (\ref{1}) with $c_{nj}=0$ for $n=1,2$ and $j>3$), and \cite%
{CP2019} \cite{CP2020} (which treat a multiplicity of analogous models);
because these papers have motivated this research and also because from them
relevant previous references can be traced. $\blacksquare $

\textbf{Remark 1-2}. The system (\ref{1}) is clearly \textit{invariant}
under the symmetry transformation
\end{subequations}
\begin{eqnarray}
&&c_{11}\Leftrightarrow c_{23}~,~c_{12}\Leftrightarrow
c_{22}~,~c_{13}\Leftrightarrow c_{21}~,~c_{14}\Leftrightarrow
c_{25}~,~c_{15}\Leftrightarrow c_{24}~,~c_{16}\Leftrightarrow c_{26}~;
\notag \\
&&x_{1}\left( t\right) \Leftrightarrow x_{2}\left( t\right) ~.
\label{1Trans}
\end{eqnarray}%
$\blacksquare $

In \textbf{Section 2} the technique used in this paper to solve the \textit{%
initial-values} problem of system (\ref{1}) is described: it involves the
introduction of $10$ parameters $A_{nm}$ and $a_{n\ell }$, in terms of which
the $12$ parameters $c_{nj}$ are expressed; the \textit{inverse} problem to
express these $10$ parameters $A_{nm}$ and $a_{n\ell }$ in terms of the $12$
\textit{a priori arbitrary} parameters $c_{nj}$ is then solved in \textbf{%
Section 3}, and it is shown there that this entails---as it were, \textit{a
posteriori}---that the $12$ parameters $c_{nj}$ are thereby required to
satisfy $4$ rather simple \textit{constraints}, which are \textit{explicitly}
exhibited. The reader who is only interested in the main results may jump
over these $2$ sections and go directly to \textbf{Section 4}, where a
\textit{summary} of the main results of this paper is presented. The
subsequent \textbf{Section 5} is devoted to two special cases of the system (%
\ref{1}) which deserve a separate treatment; and an extremely terse \textbf{%
Section 6} completes the main body of the paper; which also includes $3$
short Appendices.

\section{The technique to solve the system (\protect\ref{1})}

\textbf{First position}:
\begin{subequations}
\label{xnyn}
\begin{equation}
x_{n}\left( t\right) =A_{n1}y_{1}\left( t\right) +A_{n2}y_{2}\left( t\right)
~,  \label{xnyna}
\end{equation}%
namely%
\begin{equation}
x_{1}\left( t\right) =A_{11}y_{1}\left( t\right) +A_{12}y_{2}\left( t\right)
~,  \label{xnynb}
\end{equation}%
\begin{equation}
x_{2}\left( t\right) =A_{21}y_{1}\left( t\right) +A_{22}y_{2}\left( t\right)
~.  \label{xnync}
\end{equation}%
This assignment implies the introduction of the $4$, \textit{a priori
arbitrary}, time-independent parameters $A_{nm}$ ($n=1,2;$ $m=1,2)$ and of
the $2$ \textit{auxiliary} variables $y_{1}\left( t\right) $ and $%
y_{2}\left( t\right) $.

\textbf{Remark 2-1}. This assignment is clearly invariant under the
transformation
\end{subequations}
\begin{equation}
x_{1}\left( t\right) \Leftrightarrow x_{2}\left( t\right) ~,~~y_{1}\left(
t\right) \Leftrightarrow y_{2}\left( t\right) ~,~~A_{11}\Leftrightarrow
A_{22}~,~~A_{12}\Leftrightarrow A_{21}~.  \label{2Transx1x2}
\end{equation}%
$\blacksquare $

\textbf{Remark 2-2}. In the following (except in \textbf{Subsection 5.1}) we
generally assume that none of the $4$ parameters $A_{nm}$ vanishes, $%
A_{nm}\neq 0$. $\blacksquare $

\textbf{Remark 2-3}. The addition of $2$ \textit{a priori arbitrary}
additional parameters---say $A_{1}$ respectively $A_{2}$---to the right-hand
sides of the $2$ eqs. (\ref{xnyn}) would only complicate the following
developments without providing any significant additional generality to our
treatment. $\blacksquare $

\textbf{Evolution of the auxiliary variables }$y_{n}\left( t\right) $. Let
us hereafter assume that $y_{1}\left( t\right) $ and $y_{2}\left( t\right) $
evolve according to the following system of $2$ \textit{decoupled} ODEs:
\begin{subequations}
\label{yndot}
\begin{equation}
\dot{y}_{n}\left( t\right) =a_{n2}\left[ y_{n}\left( t\right) \right]
^{2}+a_{n1}y_{n}\left( t\right) +a_{n0}~,  \label{yndota}
\end{equation}%
namely
\begin{equation}
\dot{y}_{1}\left( t\right) =a_{12}\left[ y_{1}\left( t\right) \right]
^{2}+a_{11}y_{1}\left( t\right) +a_{10}~,  \label{yndotb}
\end{equation}%
\begin{equation}
\dot{y}_{2}\left( t\right) =a_{22}\left[ y_{2}\left( t\right) \right]
^{2}+a_{21}y_{2}\left( t\right) +a_{20}~.  \label{yndotc}
\end{equation}%
The \textit{initial-values} problem of this system of $2$ (decoupled)
ODEs---which involve the $6$ \textit{a priori arbitrary} time-independent
parameters $a_{n\ell }$ ($n=1,2;$ $\ell =2,1,0$)---is easily seen to be
\textit{explicitly solvable }(see \textbf{Appendix A}):
\end{subequations}
\begin{subequations}
\label{ynt}
\begin{equation}
y_{n}\left( t\right) =\frac{y_{n}\left( 0\right) \left[ y_{n+}-y_{n-}\exp
\left( \beta _{n}t\right) \right] -y_{n+}y_{n-}\left[ 1-\exp \left( \beta
_{n}t\right) \right] }{y_{n+}\exp \left( \beta _{n}t\right)
-y_{n-}+y_{n}\left( 0\right) \left[ 1-\exp \left( \beta _{n}t\right) \right]
}~,~~~n=1,2~,  \label{2ynt}
\end{equation}%
with the $6$ (time-independent) parameters $y_{n\pm }$ and $\beta _{n}$
defined (above and hereafter) in terms of the $6$ parameters $a_{n\ell }$ as
follows:
\begin{equation}
y_{n\pm }=\left( -a_{n1}\pm \beta _{n}\right) /\left( 2a_{n2}\right)
~,~~~\beta _{n}=\sqrt{\left( a_{n1}\right) ^{2}-4a_{n0}a_{n2}}~,~~~n=1,2~.
\label{ynplusminusbetan}
\end{equation}

\textbf{Remark 2-4}. The system of $2$ decoupled ODEs (\ref{yndot}) is
clearly invariant under the transformation
\end{subequations}
\begin{equation}
y_{1}\left( t\right) \Leftrightarrow y_{2}\left( t\right) ~,~~~a_{1\ell
}\Leftrightarrow a_{2\ell }~,~~~\ell =2,1,0~.  \label{Transy1y2a1la2l}
\end{equation}

Hence the combination of this invariance property with that reported
above---see \textbf{Remark 2-1}---clearly entails the overall invariance
property of the $2$ systems of ODEs (\ref{1}) and (\ref{yndot}), as well as
the change of variables (\ref{xnyn}), under the following transformations:%
\begin{eqnarray}
x_{1}\left( t\right) &\Leftrightarrow &x_{2}\left( t\right) ~,~~y_{1}\left(
t\right) \Leftrightarrow y_{2}\left( t\right) ~,~~A_{11}\Leftrightarrow
A_{22}~,~~A_{12}\Leftrightarrow A_{21}~,  \notag \\
a_{11} &\Leftrightarrow &a_{21}~,~~a_{12}\Leftrightarrow
a_{22}~,~~a_{10}\Leftrightarrow a_{20}~.~~~\blacksquare  \label{23Trasns}
\end{eqnarray}

Let us emphasize that, while the solutions (\ref{ynt}) are \textit{rather
simple}, their behaviors as functions of $t$---even for \textit{real} $t$
("time")---can be \textit{fairly complicated} if the parameters $\beta _{n}$
are themselves \textit{not real}; as may well be the case even if \textit{all%
} the parameters $a_{n\ell }$ are \textit{real }numbers: see (\ref%
{ynplusminusbetan}). On the other hand if both $\beta _{n}=\mathbf{i}\rho
_{n}\omega $ with $n=1,2$, $\omega $ an \textit{arbitrary (nonvanishing) real%
} number and with \textit{both} parameters $\rho _{n}$ \textit{rational
(nonvanishing) real} numbers, then clearly (or, if need be, see for instance
\cite{C2008}) the evolution of both the pairs $y_{n}\left( t\right) $ and $%
x_{n}\left( t\right) $---as functions of \textit{real} $t$ ("time")---is
\textit{completely periodic} with a period \textit{independent }of the
respective initial data, which is an \textit{integer }multiple of the basic
period $2\pi /\left\vert \omega \right\vert $ ("isochrony").

Let us also note that $y_{n\pm }$ (see (\ref{ynplusminusbetan})) are the
\textit{equilibrium} positions of the system (\ref{yndot}), while the
\textit{asymptotic} behavior of $y_{n}\left( t\right) $ as the (\textit{real}%
) time $t$ \textit{diverges} is always rather simple: indeed if $Re%
\left[ \beta _{n}\right] <0$, then
\begin{subequations}
\begin{equation}
\underset{t\rightarrow +\infty }{\lim }\left[ y_{n}\left( t\right) \right]
=y_{n+}~;  \label{yn+}
\end{equation}%
while if $Re\left[ \beta _{n}\right] >0$ then
\begin{equation}
\underset{t\rightarrow +\infty }{\lim }\left[ y_{n}\left( t\right) \right]
=y_{n-}~.  \label{yn-}
\end{equation}%

And clearly the corresponding \textit{equilibrium} configurations and
\textit{asymptotic} behaviors of the variables $x_{n}\left( t\right) $ are
as well rather simple (see (\ref{xnyn})): but note that the asymptotic
behavior of these variables $x_{n}\left( t\right) $ is \textit{%
asymptotically isochronous} (see \cite{CGU2008}), when one and only one of
the $2$ quantities $\beta _{n}$ is \textit{purely imaginary}.

It is of course evident that the \textit{coupled} system of $2$ ODEs implied
by the relations (\ref{xnyn}) and by the $2$ ODEs (\ref{yndot}) satisfied by
the $2$ variables $y_{n}\left( t\right) $ is \textit{identical} to the
system of $2$ ODEs (\ref{1}) satisfied by the $2$ variables $x_{n}\left(
t\right) $, of course provided the $12$ parameters $c_{nj}$ are \textit{%
appropriately} expressed in terms of the $10=4+6$ parameters $A_{nm}$ in (%
\ref{xnyn}) and $a_{n\ell }$ in (\ref{yndot}). The corresponding computation
of the explicit expressions of the $12$ parameters $c_{nj}$ in terms of the $%
10$ parameters $A_{nm}$, $a_{n\ell }$ is a standard---if tedious---task (see
\textbf{Appendix B}), yielding (for the $6$ parameters $c_{1j}$) the
following results:
\end{subequations}
\begin{subequations}
\label{c1j}
\begin{equation}
c_{11}=\left[ a_{12}A_{11}\left( A_{22}\right) ^{2}+a_{22}A_{12}\left(
A_{21}\right) ^{2}\right] /D^{2}~,
\end{equation}%
\begin{equation}
c_{12}=-2A_{11}A_{12}\left[ a_{12}A_{22}+a_{22}A_{21}\right] /D^{2}~,
\end{equation}%
\begin{equation}
c_{13}=A_{11}A_{12}\left[ a_{12}A_{12}+a_{22}A_{11}\right] /D^{2}~,
\end{equation}%
\begin{equation}
c_{14}=\left( a_{11}A_{11}A_{22}-a_{21}A_{12}A_{21}\right) /D~,
\end{equation}%
\begin{equation}
c_{15}=-\left( a_{11}-a_{21}\right) A_{11}A_{12}/D~,
\end{equation}%
\begin{equation}
c_{16}=a_{10}A_{11}+a_{20}A_{12}~,
\end{equation}%
where the quantity $D$ is defined, above and hereafter, as follows:
\end{subequations}
\begin{equation}
D=A_{11}A_{22}-A_{12}A_{21}~.  \label{AA}
\end{equation}

The analogous formulas for the $6$ parameters $c_{2j}$ can be obtained from
those written above via the transformations (see \textbf{Remarks 1-2},
\textbf{2-1 }and\textbf{\ 2-2})\textbf{\ }
\begin{eqnarray}
c_{11} &\Leftrightarrow &c_{23}~,~c_{12}\Leftrightarrow
c_{22}~,~c_{13}\Leftrightarrow c_{21}~,~c_{14}\Leftrightarrow
c_{25}~,~c_{15}\Leftrightarrow c_{24}~,~c_{16}\Leftrightarrow c_{26}~;
\notag \\
A_{11} &\Leftrightarrow &A_{22}~,~~~A_{12}\Leftrightarrow
A_{21}~;~~~a_{1\ell }\Leftrightarrow a_{2\ell }~,  \label{2TransfcnjAnmanl}
\end{eqnarray}%
under which the quantity $D$, see (\ref{AA}), is clearly invariant.

Hence they read as follows:
\begin{subequations}
\label{c2j}
\begin{equation}
c_{23}=\left[ a_{22}A_{22}\left( A_{11}\right) ^{2}+a_{12}A_{21}\left(
A_{12}\right) ^{2}\right] /D^{2}~,
\end{equation}%
\begin{equation}
c_{22}=-2A_{22}A_{21}\left[ a_{22}A_{11}+a_{12}A_{12}\right] /D^{2}~,
\end{equation}%
\begin{equation}
c_{21}=A_{22}A_{21}\left[ a_{22}A_{21}+a_{12}A_{22}\right] /D^{2}~,
\end{equation}%
\begin{equation}
c_{25}=\left( a_{21}A_{11}A_{22}-a_{11}A_{12}A_{21}\right) /D~,
\end{equation}%
\begin{equation}
c_{24}=-\left( a_{21}-a_{11}\right) A_{22}A_{21}/D~,
\end{equation}%
\begin{equation}
c_{26}=a_{20}A_{22}+a_{10}A_{21}~.
\end{equation}

These findings of course imply that the solution of the \textit{%
initial-values} problem for the system (\ref{1}) is provided via the
formulas (\ref{xnyn}) from the \textit{explicit} solutions (\ref{ynt}) of
the \textit{initial-values} problem for the system (\ref{yndot}), hence
essentially via rather simple, quite \textit{explicit,} \textit{algebraic}
operations; this can be done for any assignment of the $12$ parameters $%
c_{nj}$, such that the \textit{explicit} formulas (\ref{c1j}) and (\ref{c2j}%
)---expressing the $12$ parameters $c_{nj}$ in terms of the $10$ parameters $%
A_{nm}$ and $a_{n\ell }$---can be \textit{inverted}: in the following
\textbf{Section 3} we show how this can be done, provided the $12$
parameters $c_{nj}$ satisfy $4$ \textit{constraints}.

\section{The inverse problem: expressing the $10$ parameters $A_{nm}$ and $%
a_{n\ell }$ in terms of the $12$ coefficients $c_{nj}$}

Our main task in this \textbf{Section 3} is to discuss the \textit{inversion}
of the $12$ relations obtained above---see (\ref{c1j}) and (\ref{c2j})%
\textbf{\ }---expressing the $12$ coefficients $c_{nj}$ in terms of the $10$
parameters $A_{nm}$ and $a_{n\ell }$; namely to show how, given $12$ \textit{%
a priori arbitrary }coefficients $c_{nj}$, the $10$ corresponding parameters
$A_{nm}$ and $a_{n\ell }$ can be computed. We show below how this task can
be \textit{explicitly }accomplished; but that it entails that the $12$
coefficients $c_{nj}$ ($n=1,2$; $j=1,2,3,4,5,6$) must---as it were, \textit{%
a posteriori}---satisfy $4$ \textit{algebraic} conditions (\textit{explicitly%
} obtained below).

An obvious route to achieve our main task is to try and solve the $12$
\textit{algebraic} equations (\ref{c1j}) and (\ref{c2j}) for the $10$
parameters $A_{nm}$ and $a_{n\ell }$; but given the fairly \textit{large}
number of these \textit{algebraic} relations and their \textit{nonlinear}
character this is a nontrivial job (beyond the power of standard algebraic
manipulation packages such as \textbf{Mathematica}, used on a modern PC).
More progress in this direction can be made via the following alternative
procedure.

Let us note that the relation (\ref{y1xn}) in \textbf{Appendix B} implies
that the variable $y_{1}\left( t\right) $ clearly satisfies the ODE
\end{subequations}
\begin{subequations}
\begin{equation}
\dot{y}_{1}=\left( A_{22}\dot{x}_{1}-A_{12}\dot{x}_{2}\right) /D~,
\end{equation}%
hence, via (\ref{1}),%
\begin{eqnarray}
&&\dot{y}_{1}=\left\{ A_{22}\left[ c_{11}\left( x_{1}\right)
^{2}+c_{12}x_{1}x_{2}+c_{13}\left( x_{2}\right)
^{2}+c_{14}x_{1}+c_{15}x_{2}+c_{16}\right] \right.  \nonumber \\
&&\left. -A_{12}\left[ c_{21}\left( x_{1}\right)
^{2}+c_{22}x_{1}x_{2}+c_{23}\left( x_{2}\right)
^{2}+c_{24}x_{1}+c_{25}x_{2}+c_{26}\right] \right\} /D~,  \nonumber \\
&&
\end{eqnarray}%
yielding, via (\ref{xnyn}) and some trivial if tedious algebra,
\end{subequations}
\begin{equation}
\dot{y}_{1}=b_{11}\left( y_{1}\right) ^{2}+b_{12}y_{1}y_{2}+b_{13}\left(
y_{2}\right) ^{2}+b_{14}y_{1}+b_{15}y_{2}+b_{16}~,  \label{3y1dot}
\end{equation}%
with
\begin{subequations}
\label{3b1j}
\begin{eqnarray}
b_{11} &=&\left[ \left( A_{11}\right) ^{2}\left(
-A_{12}c_{21}+A_{22}c_{11}\right) +A_{11}A_{21}\left(
-A_{12}c_{22}+A_{22}c_{12}\right) \right.  \nonumber \\
&&\left. +\left( A_{21}\right) ^{2}\left( -A_{12}c_{23}+A_{22}c_{13}\right)
\right] /D~,
\end{eqnarray}%
\begin{eqnarray}
&&b_{12}=\left\{ A_{11}\left[ -2\left( A_{12}\right)
^{2}c_{21}+A_{12}A_{22}\left( 2c_{11}-c_{22}\right) +\left( A_{22}\right)
^{2}c_{12}\right] \right.  \nonumber \\
&&\left. +A_{21}\left[ -\left( A_{12}\right) ^{2}c_{22}+A_{12}A_{22}\left(
c_{12}-2c_{23}\right) +2\left( A_{22}\right) ^{2}c_{13}\right] \right\} /D~,
\end{eqnarray}%
\begin{eqnarray}
&&b_{13}=\left\{ \left( A_{12}\right) ^{2}\left[ -A_{12}c_{21}+A_{22}\left(
c_{11}-c_{22}\right) \right] \right.  \nonumber \\
&&\left. +\left( A_{22}\right) ^{2}\left[ A_{12}\left( c_{12}-c_{23}\right)
+A_{22}c_{13}\right] \right\} /D~,
\end{eqnarray}%
\begin{equation}
b_{14}=\left[ A_{11}\left( -A_{12}c_{24}+A_{22}c_{14}\right) +A_{21}\left(
-A_{12}c_{25}+A_{22}c_{15}\right) \right] /D~,  \nonumber
\end{equation}

\begin{equation}
b_{15}=\left[ -\left( A_{12}\right) ^{2}c_{24}+A_{12}A_{22}\left(
c_{14}-c_{25}\right) +\left( A_{22}\right) ^{2}c_{15}\right] /D~,
\end{equation}

\begin{equation}
b_{16}=\left( -A_{12}c_{26}+A_{22}c_{16}\right) /D~.
\end{equation}

And now a comparison of (\ref{yndotb}) with (\ref{3y1dot}) yields, via (\ref%
{3b1j}), the following $6$ relations:
\end{subequations}
\begin{subequations}
\label{3a1L}
\begin{eqnarray}
a_{12} &=&\left[ \left( A_{11}\right) ^{2}\left(
-A_{12}c_{21}+A_{22}c_{11}\right) +A_{11}A_{21}\left(
-A_{12}c_{22}+A_{22}c_{12}\right) \right.  \nonumber \\
&&\left. +\left( A_{21}\right) ^{2}\left( -A_{12}c_{23}+A_{22}c_{13}\right)
\right] /D~,
\end{eqnarray}

\begin{equation}
a_{11}=\left[ A_{11}\left( -A_{12}c_{24}+A_{22}c_{14}\right) +A_{21}\left(
-A_{12}c_{25}+A_{22}c_{15}\right) \right] /D~,
\end{equation}%
\begin{equation}
a_{10}=\left( -A_{12}c_{26}+A_{22}c_{16}\right) /D~;
\end{equation}%
\end{subequations}
\begin{subequations}
\label{3Anmc1L}
\begin{eqnarray}
&&A_{11}\left[ -2\left( A_{12}\right) ^{2}c_{21}+A_{12}A_{22}\left(
2c_{11}-c_{22}\right) +\left( A_{22}\right) ^{2}c_{12}\right]  \nonumber \\
&&+A_{21}\left[ -\left( A_{12}\right) ^{2}c_{22}+A_{12}A_{22}\left(
c_{12}-2c_{23}\right) +2\left( A_{22}\right) ^{2}c_{13}\right] =0~,
\label{3Anm1}
\end{eqnarray}%
\begin{eqnarray}
&&\left( A_{12}\right) ^{2}\left[ -A_{12}c_{21}+A_{22}\left(
c_{11}-c_{22}\right) \right]  \nonumber \\
&&+\left( A_{22}\right) ^{2}\left[ A_{12}\left( c_{12}-c_{23}\right)
+A_{22}c_{13}\right] =0~,  \label{3Anm2}
\end{eqnarray}%
\begin{equation}
-\left( A_{12}\right) ^{2}c_{24}+A_{12}A_{22}\left( c_{14}-c_{25}\right)
+\left( A_{22}\right) ^{2}c_{15}=0~.  \label{FirstEqA1A2}
\end{equation}

By a completely analogous development, based on the ODE (\ref{yndotc})
satisfied by $y_{2}\left( t\right) $ rather than (\ref{yndotb}) satisfied by
$y_{1}\left( t\right) $---or, more easily, via the symmetry properties
associated to the transformation (\ref{2TransfcnjAnmanl})---one gets the
following $6$ additional relations:
\end{subequations}
\begin{subequations}
\label{3a2L}
\begin{eqnarray}
a_{22} &=&\left[ \left( A_{22}\right) ^{2}\left(
-A_{21}c_{13}+A_{11}c_{23}\right) +A_{22}A_{12}\left(
-A_{21}c_{12}+A_{11}c_{22}\right) \right.  \nonumber \\
&&\left. +\left( A_{12}\right) ^{2}\left( -A_{21}c_{11}+A_{11}c_{21}\right)
\right] /D~,
\end{eqnarray}%
\begin{equation}
a_{21}=\left[ A_{22}\left( -A_{21}c_{15}+A_{11}c_{25}\right) +A_{12}\left(
-A_{21}c_{14}+A_{11}c_{24}\right) \right] /D~,
\end{equation}%
\begin{equation}
a_{20}=\left( -A_{21}c_{16}+A_{11}c_{26}\right) /D~;
\end{equation}%
\end{subequations}
\begin{subequations}
\label{Anm4}
\begin{eqnarray}
&&A_{22}\left[ -2\left( A_{21}\right) ^{2}c_{13}+A_{21}A_{11}\left(
2c_{23}-c_{12}\right) +\left( A_{11}\right) ^{2}c_{22}\right]  \nonumber \\
&&+A_{12}\left[ -\left( A_{21}\right) ^{2}c_{12}+A_{21}A_{11}\left(
c_{22}-2c_{11}\right) +2\left( A_{11}\right) ^{2}c_{21}\right] =0~,
\label{3Anm3}
\end{eqnarray}%
\begin{eqnarray}
&&\left( A_{21}\right) ^{2}\left[ -A_{21}c_{13}+A_{11}\left(
c_{23}-c_{12}\right) \right]  \nonumber \\
&&+\left( A_{11}\right) ^{2}\left[ A_{21}\left( c_{22}-c_{11}\right)
+A_{11}c_{21}\right] =0~,  \label{3Anm4}
\end{eqnarray}%
\begin{equation}
-\left( A_{21}\right) ^{2}c_{15}+A_{21}A_{11}\left( c_{25}-c_{14}\right)
+\left( A_{11}\right) ^{2}c_{24}=0~.  \label{SecondEqA1A2}
\end{equation}

It is thus seen that the $6$ parameters $a_{n\ell }$ ($n=1,2;$ $\ell =2,1,0)$
are given \textit{explicitly} by the $6$ formulas (\ref{3a1L}) and (\ref%
{3a2L}) in terms of the $12$ coefficients $c_{nj}$ and the $4$ parameters $%
A_{nm}$.

The remaining task is to extract the expressions of the $4$ parameters $%
A_{nm}$ in terms of the $6$ parameters $c_{nk}$ ($n=1,2;$ $k=1,2,3$), the
only ones featured in the remaining $6$ \textit{algebraic} equations (\ref%
{3Anmc1L}) and (\ref{Anm4}).

To this end, let us now introduce the $2$ \textit{auxiliary} parameters $%
z_{1}$ and $z_{2}$:
\end{subequations}
\begin{equation}
z_{1}=A_{11}/A_{21}~,~~~z_{2}=A_{12}/A_{22}~;  \label{3z1z2}
\end{equation}%
it is then easily seen that the $2$ eqs. (\ref{FirstEqA1A2}) and (\ref%
{SecondEqA1A2}) yield---recall \textbf{Remark 2-2}---the \textit{same cubic }%
equation for these $2$ quantities:
\begin{subequations}
\label{3EqzCubic}
\begin{equation}
c_{21}\left( z_{n}\right) ^{3}+\left( c_{22}-c_{11}\right) \left(
z_{n}\right) ^{2}+\left( c_{23}-c_{12}\right) z_{n}-c_{13}=0~,~~~n=1,2~,
\label{3Eqz}
\end{equation}%
namely%
\begin{equation}
c_{21}\left( z_{1}\right) ^{3}+\left( c_{22}-c_{11}\right) \left(
z_{1}\right) ^{2}+\left( c_{23}-c_{12}\right) z_{1}-c_{13}=0~,  \label{3Eqz1}
\end{equation}%
\begin{equation}
c_{21}\left( z_{2}\right) ^{3}+\left( c_{22}-c_{11}\right) \left(
z_{2}\right) ^{2}+\left( c_{23}-c_{12}\right) z_{2}-c_{13}=0~.  \label{3Eqz2}
\end{equation}

\textbf{Remark 3-1}. Note that, consistently with the transformations (\ref%
{1Trans}) and (\ref{23Trasns}), the corresponding transformations of the $2$
auxiliary parameters $z_{n}$ are $z_{1}\Leftrightarrow 1/z_{2}$ and (of
course) $z_{2}\Leftrightarrow 1/z_{1}$, implying the invariance under all
these transformations of the eqs. (\ref{3EqzCubic}). $\blacksquare $

The $2$ algebraic equations (\ref{3EqzCubic}) allow to compute (explicitly,
via the Cardano formulas) the $2$ quantities $z_{n}$ in terms of the $6$
coefficients $c_{nk}$; of course, they do \textit{not} imply that the $2$
quantities $z_{1}$ and $z_{2}$ coincide, indeed we \textit{exclude}
hereafter this possibility because it would imply the vanishing of $D$ (see (%
\ref{AA}) and (\ref{3z1z2})).

But more progress is possible.

Indeed, let us take advantage of the definitions (\ref{3z1z2}) to rewrite
the $2$ eqs. (\ref{3Anm1}) and (\ref{3Anm3}), getting thereby (again
recalling \textbf{Remark 2-2})
\end{subequations}
\begin{subequations}
\label{3Eqszc}
\begin{eqnarray}
&&\left[ -2c_{21}\left( z_{n+1}\right) ^{2}+\left( 2c_{11}-c_{22}\right)
z_{n+1}+c_{12}\right] z_{n}  \nonumber \\
&&-c_{22}\left( z_{n+1}\right) ^{2}+\left( c_{12}-2c_{23}\right)
z_{n+1}+2c_{13}=0~,~~n=1,2~\mod \left[ 2\right] ~,
\end{eqnarray}%
namely (also dividing by $2$)%
\begin{eqnarray}
&&\left[ -c_{21}\left( z_{2}\right) ^{2}+\left( c_{11}-c_{22}/2\right)
z_{2}+c_{12}/2\right] z_{1}  \nonumber \\
&&-c_{22}\left( z_{2}\right) ^{2}/2+\left( c_{12}/2-c_{23}\right)
z_{2}+c_{13}=0~,  \label{3FirstConstraint}
\end{eqnarray}%
\begin{eqnarray}
&&\left[ -c_{21}\left( z_{1}\right) ^{2}+\left( c_{11}-c_{22}/2\right)
z_{1}+c_{12}/2\right] z_{2}  \nonumber \\
&&-c_{22}\left( z_{1}\right) ^{2}/2+\left( c_{12}/2-c_{23}\right)
z_{1}+c_{13}=0~.  \label{3SecondConstraint}
\end{eqnarray}

Let us now sum the $2$ eqs. (\ref{3Eqz1}) and (\ref{3SecondConstraint}), and
likewise the $2$ eqs. (\ref{3Eqz2}) and (\ref{3FirstConstraint}). We thus
obtain (using the fact that $z_{1}-z_{2}\neq 0$; see above) the same \textit{%
quadratic} equation for the $2$ quantities $z_{1}$ and $z_{2}$:
\end{subequations}
\begin{subequations}
\begin{equation}
2c_{21}\left( z_{n}\right) ^{2}-\left( 2c_{11}-c_{22}\right)
z_{n}-c_{12}=0~,~~~n=1,2~,  \label{3Eqzn}
\end{equation}%
implying
\begin{equation}
z_{n}=\left[ 2c_{11}-c_{22}+\left( -1\right) ^{n}\sqrt{\left(
2c_{11}-c_{22}\right) ^{2}+8c_{12}c_{21}}\right] /\left( 4c_{21}\right)
~,~~~n=1,2~.  \label{3zn}
\end{equation}%
These formulas (\ref{3zn}) feature of course only \textit{square-roots}%
---rather than the \textit{cubic-roots} that would be featured by the
Cardano solutions of the \textit{cubic} equations (\ref{3EqzCubic})---and
moreover they yield the \textit{explicit} expressions (\ref{3zn}) of the $2$
auxiliary parameters $z_{n}$ in terms of (only!) the $4$ parameters $c_{nm}$
($n=1,2$; $m=1,2$). Hence by inserting these expressions of $z_{1}$ and $%
z_{2}$ in any $2$ of the $4$ eqs. (\ref{3EqzCubic}) and (\ref{3Eqszc}), we
get a system of $2$ \textit{algebraic} equations satisfied by the $6$
parameters $c_{nk}$ ($n=1,2$; $k=1,2,3$) which features the $2$ parameters $%
c_{13}$ and $c_{23}$ (only!) \textit{linearly} and therefore allows to
express both these coefficients \textit{explicitly} in terms of the other $4$
coefficients $c_{nm}$ ($n=1,2;m=1,2$). For instance the $2$ eqs. (\ref%
{3FirstConstraint}) and (\ref{3SecondConstraint}) yield (by subtracting the
second multiplied by $z_{1}$ from the first multiplied by $z_{2},$ and by
subtracting the second from the first)
\end{subequations}
\begin{subequations}
\begin{equation}
c_{13}=-c_{11}z_{1}z_{2}-c_{12}\left( z_{1}+z_{2}\right) /2~,
\end{equation}%
\begin{equation}
c_{23}=-c_{21}z_{1}z_{2}-c_{22}\left( z_{1}+z_{2}\right) /2~,
\end{equation}%
hence, via (\ref{3zn}),
\end{subequations}
\begin{subequations}
\label{3ConstraintsFinal}
\begin{equation}
4c_{13}c_{21}-c_{12}c_{22}=0~,
\end{equation}%
\begin{equation}
2\left( -c_{12}+2c_{23}\right) c_{21}+\left( 2c_{11}-c_{22}\right) c_{22}=0~.
\end{equation}

\textbf{Remark 3-2}. Since throughout our treatment we have assumed that $%
z_{1}$ is different from $z_{2},$ clearly the formula (\ref{3zn}) implies
that the $4$ parameters $c_{11},$ $c_{12},$ $c_{21},$ $c_{22}$ must satisfy
the \textit{inequality}
\end{subequations}
\begin{equation}
\left( 2c_{11}-c_{22}\right) ^{2}+8c_{12}c_{21}\neq 0~.  \label{3Ineq1}
\end{equation}%
$\blacksquare $

Two additional relations can be obtained by inserting in the $2$ eqs. (\ref%
{FirstEqA1A2}) and (\ref{SecondEqA1A2}) the expressions%
\begin{equation}
A_{11}=z_{1}A_{21}~,~~~A_{12}=z_{2}A_{22}~,
\end{equation}%
implied by (\ref{3z1z2}), obtaining thereby (again) $2$ identical \textit{%
second-degree} equations for the $2$ parameters $z_{1}$ and $z_{2}:$%
\begin{subequations}
\label{3cc1245}
\begin{equation}
c_{24}\left( z_{n}\right) ^{2}+\left( c_{25}-c_{14}\right)
z_{n}-c_{15}=0~,~~~n=1,2~,  \label{3SecondEqzn}
\end{equation}%
implying of course%
\begin{equation}
z_{n}=\left[ c_{14}-c_{25}+\left( -1\right) ^{n}\sqrt{\left(
c_{14}-c_{25}\right) ^{2}+4c_{15}c_{24}}\right] /\left( 2c_{24}\right)
~,~~~n=1,2~.  \label{3znSecond}
\end{equation}

\textbf{Remark 3-3}. Note that these equations (\ref{3cc1245}) only involve
the $4$ parameters $c_{14},$ $c_{24},$ $c_{15},$ $c_{25};$ and that the
condition $z_{1}\neq z_{2}$ entails the \textit{inequality}
\end{subequations}
\begin{equation}
\left( c_{25}-c_{14}\right) ^{2}+4c_{15}c_{24}\neq 0~.  \label{3Ineq2}
\end{equation}%
$\blacksquare $

Since from the $4$ eqs. (\ref{3Anm1}), (\ref{3Anm2}), (\ref{3Anm3}), (\ref%
{3Anm4}) we extracted the $2$ \textit{constraints }(\ref{3ConstraintsFinal}%
), clearly from these $4$ equations we can only obtain $2$ additional
relations constraining the parameters $c_{nj}$. A convenient way to get such
relations is to subtract the eq. (\ref{3SecondEqzn}) multiplied by $2c_{21}$
from the eq. (\ref{3Eqzn}) itself multiplied by $c_{24},$ getting thereby
the following $2$ identical \textit{first-degree} equations for the
parameters $z_{1}$ and $z_{2}$:%
\begin{eqnarray}
&&-\left[ c_{24}\left( 2c_{11}-c_{22}\right) +2c_{21}\left(
c_{25}-c_{14}\right) \right] z_{n}  \nonumber \\
&=&c_{12}c_{24}-2c_{15}c_{21}~,~~~n=1,2~.  \label{3FirstDegreeEq}
\end{eqnarray}%
But these $2$ first-degree equations seem to imply that $z_{1}=z_{2},$ while
we know that this is \textit{not} the case (at least, provided the two
inequalities (\ref{3Ineq1}) and (\ref{3Ineq2}) hold true; which is generally
the case for any \textit{generic} assignment of $12$ parameters $c_{nj}$).
Hence these \textit{first-degree} eqs. (\ref{3FirstDegreeEq}) satisfied by $%
z_{1}$ and $z_{2}$ must have the property to be \textit{identically}
satisfied for any arbitrary value of $z_{1}$ and $z_{2},$ which is of course
the case provided the $8$ parameters $c_{11},$ $c_{12},$ $c_{14},$ $c_{15},$
$c_{21},$ $c_{22},$ $c_{24},$ $c_{25}$ satisfy the following $2$ \textit{%
constraints}:
\begin{subequations}
\label{3ThirdFourthCon}
\begin{equation}
c_{24}\left( 2c_{11}-c_{22}\right) +2c_{21}\left( c_{25}-c_{14}\right) =0~,
\end{equation}%
\begin{equation}
c_{12}c_{24}-2c_{15}c_{21}=0
\end{equation}%
(to obtain the first of these $2$ equations we assumed $c_{24}\neq 0$,
consistently with our assumption that the $12$ parameters $c_{nj}$ have
\textit{generic} values).

In conclusion, we have obtained $4$ \textit{constraints} on the $10$
parameters $c_{np}$ ($n=1,2$; $p=1,2,3,4,5$): see the $2$ eqs. (\ref%
{3ConstraintsFinal}) and the $2$ eqs. (\ref{3ThirdFourthCon}); note that the
$2$ parameters $c_{n6}$ are \textit{not} involved at all in these \textit{%
constraints}.

There remain to compute the $4$ parameters $A_{nm}$.

To compute the $4$ parameters $A_{nm}$, rather than using the $6$ eqs. (\ref%
{3Anmc1L}) and (\ref{Anm4})---out of which we already extracted the $4$
\textit{constraints} (\ref{3ConstraintsFinal}) and (\ref{3ThirdFourthCon});
so that we can expect to be only able to compute only $2$ of the $4$
parameters $A_{nm}$ in terms of the other $2$ (and of course the
coefficients $c_{nj}$)---the simplest way is to use the relations implied by
the definitions (\ref{3z1z2}):
\end{subequations}
\begin{equation}
A_{11}=z_{1}A_{21}~,~~~A_{12}=z_{2}A_{22}~;  \label{3AA}
\end{equation}%
here of course $z_{1}$ and $z_{2}$ are defined by their expressions (\ref%
{3zn}) or, equivalently, (\ref{3znSecond}), and the $2$ parameters $A_{21}$
and $A_{22}$ can be considered as \textit{free} parameters; so that these
relations can be rewritten as follows:%
\begin{equation}
A_{21}=\lambda _{1}~,~~~A_{22}=\lambda _{2}~,~~~A_{11}=z_{1}\lambda
_{1}~,~~~A_{12}=z_{2}\lambda _{2}~,  \label{3Alanda}
\end{equation}%
with $\lambda _{1}$ and $\lambda _{2}$ two \textit{arbitrary} (nonvanishing)
parameters.

This concludes both the identification of the \textit{subclass} of the
dynamical system (\ref{1}) which is treated in this paper and the solution
of its \textit{initial-values} problem; except for the further step of
inserting in all the relevant formulas---in addition to the expression (\ref%
{3Alanda})---the following rather simple expressions, say, of $c_{13}$ and $%
c_{23},$%
\begin{subequations}
\begin{equation}
c_{13}=c_{12}c_{22}/\left( 4c_{21}\right) ~,
\end{equation}%
\begin{equation}
c_{23}=\left[ 2c_{12}c_{21}-c_{22}\left( 2c_{11}-c_{22}\right) \right]
/\left( 4c_{21}\right)
\end{equation}%
implied by (\ref{3ConstraintsFinal}), and likewise, say, of $c_{24}$ and $%
c_{25},$%
\end{subequations}
\begin{subequations}
\begin{equation}
c_{24}=2c_{15}c_{21}/c_{12}~,
\end{equation}%
\begin{equation}
c_{25}=c_{14}-c_{24}\left( 2c_{11}-c_{22}\right) /\left( 2c_{21}\right) ~,
\end{equation}%
implied by (\ref{3ThirdFourthCon}) (we ignore the \textit{nongeneric} cases
with $c_{21}=c_{24}=0$).

An essential compendium of all the relevant formulas is displayed in the
following \textbf{Section 4}, for the convenience of the reader who is more
interested in using these findings than in following their derivation.

\textbf{Remark 3-4}. A final observation. The reader who has followed our
derivation up to this point might justifiably be puzzled by the fact that
our solution seems to feature the $2$ \textit{free} parameters $\lambda _{1}$
and $\lambda _{2}$. But in fact these $2$ \textit{free} parameters are
\textit{not} present at all in the solution $x_{n}\left( t\right) $ ($n=1,2$%
) of the dynamical system (\ref{1}) obtained above. This is proven in
\textbf{Appendix C}; the development reported there are also useful to get
the final formulas reported in the following \textbf{Section 4} (which
indeed do \textit{not} feature the $2$ parameters $\lambda _{n}$). $%
\blacksquare $

\section{A summary of the solution of the system (\protect\ref{1})}

In this \textbf{Section 4} we summarize the main results obtained in this
paper so far. For the convenience of the reader who is only interested in
these results and not in their derivation, we report these findings in a
self-consistent fashion, even at the cost of the repetition of some key
formulas already displayed above and in the \textbf{Appendices}.

Let us recall that our focus is on the system of $2$ nonlinearly-coupled
first-order ODEs (\ref{1}), i. e.
\end{subequations}
\begin{equation}
\dot{x}_{n}=c_{n1}\left( x_{1}\right) ^{2}+c_{n2}x_{1}x_{2}+c_{n3}\left(
x_{2}\right) ^{2}+c_{n4}x_{1}+c_{n5}x_{2}+c_{n6}~,~~~n=1,2~.  \label{41}
\end{equation}

Our main finding is the solution in \textit{explicit} form of the \textit{%
initial-values} problem for this system, which is however achieved only
provided its $12$ \textit{a priori arbitrary }parameters $c_{nj}$ ($n=1,2;$ $%
j=1,2,3,4,5,6$) do satisfy---as it were, \textit{a posteriori}---the
following $4$ \textit{algebraic constraints}:
\begin{subequations}
\label{4ConstraintsFinal}
\begin{equation}
4c_{13}c_{21}-c_{12}c_{22}=0~,
\end{equation}%
\begin{equation}
2\left( -c_{12}+2c_{23}\right) c_{21}+\left( 2c_{11}-c_{22}\right) c_{22}=0~,
\end{equation}%
\begin{equation}
c_{24}\left( 2c_{11}-c_{22}\right) +2c_{21}\left( c_{25}-c_{14}\right) =0~,
\end{equation}%
\begin{equation}
c_{12}c_{24}-2c_{15}c_{21}=0~.
\end{equation}%
Note that the first $2$ of these $4$ \textit{constraints} only involve the $%
6 $ parameters $c_{nk}$ ($n=1,2;$ $k=1,2,3$), and the last $2$ only involve
the $8$ parameters $c_{11},$ $c_{12},$ $c_{14},$ $c_{15},$ $c_{21},$ $%
c_{22}, $ $c_{24},$ $c_{25};$ while the $2$ parameters $c_{n6}$ are \textit{%
unconstrained} and only influence (see below) the $2$ parameters $\alpha
_{n0}$.

The \textit{explicit} solution of the \textit{initial-values} problem for
the system (\ref{41}) with (\ref{4ConstraintsFinal}) is then provided by the
following formulas, for whose derivation the interested reader should go
through the developments reported in the rest of this paper (for some
guidance see below \textbf{Remark 4-1}):
\end{subequations}
\begin{subequations}
\label{4xnt}
\begin{equation}
x_{1}\left( t\right) =z_{1}w_{1}\left( t\right) +z_{2}w_{2}\left( t\right) ~,
\label{4xnta}
\end{equation}%
\begin{equation}
x_{2}\left( t\right) =w_{1}\left( t\right) +w_{2}\left( t\right) ~;
\label{4xntb}
\end{equation}%
\end{subequations}
\begin{equation}
z_{n}=\left[ 2c_{11}-c_{22}+\left( -1\right) ^{n}\sqrt{\left(
2c_{11}-c_{22}\right) ^{2}+8c_{12}c_{21}}\right] /\left( 4c_{21}\right)
~,~~~n=1,2~;  \label{4zn}
\end{equation}%
\begin{eqnarray}
&&w_{n}\left( t\right) \equiv w_{n}\left( C,t\right) ~,  \nonumber \\
&=&\frac{w_{n}\left( 0\right) \left[ w_{n+}-w_{n-}\exp \left( \beta
_{n}t\right) \right] -w_{n+}w_{n-}\left[ 1-\exp \left( \beta _{n}t\right) %
\right] }{w_{n+}\exp \left( \beta _{n}t\right) -w_{n-}+w_{n}\left( 0\right) %
\left[ 1-\exp \left( \beta _{n}t\right) \right] }~,  \nonumber \\
n &=&1,2~;  \label{4wnt}
\end{eqnarray}%
\begin{subequations}
\label{4w1w20}
\begin{equation}
w_{1}\left( 0\right) \equiv w_{1}\left( C,0\right) =\left[ x_{1}\left(
0\right) -z_{2}x_{2}\left( 0\right) \right] /\left( z_{1}-z_{2}\right) ~,
\label{4w10}
\end{equation}%
\begin{equation}
w_{2}\left( 0\right) \equiv w_{2}\left( C,0\right) =-\left[ x_{1}\left(
0\right) -z_{1}x_{2}\left( 0\right) \right] /\left( z_{1}-z_{2}\right) ~;
\label{4w20}
\end{equation}%
\end{subequations}
\begin{equation}
w_{n\pm }=\left( -\alpha _{n1}\pm \beta _{n}\right) /\left( 2\alpha
_{n2}\right) ~,~~~\beta _{n}=\sqrt{\left( \alpha _{n1}\right) ^{2}-4\alpha
_{n0}\alpha _{n2}}~,~~n=1,2~;  \label{4wnpm}
\end{equation}%
\begin{subequations}
\label{4alphanl}
\begin{equation}
\alpha _{12}=\left[ \left( z_{1}\right) ^{2}\left( c_{11}-z_{2}c_{21}\right)
+z_{1}\left( c_{12}-z_{2}c_{22}\right) +c_{13}-z_{2}c_{23}\right] /\left(
z_{1}-z_{2}\right) ~,
\end{equation}%
\begin{equation}
\alpha _{11}=\left[ z_{1}\left( c_{14}-z_{2}c_{24}\right) +c_{15}-z_{2}c_{25}%
\right] /\left( z_{1}-z_{2}\right) ~,
\end{equation}%
\begin{equation}
\alpha _{10}=\left( c_{16}-z_{2}c_{26}\right) /\left( z_{1}-z_{2}\right) ~;
\end{equation}

\begin{eqnarray}
\alpha _{22} &=&\left[ \left( -c_{13}+z_{1}c_{23}\right) +z_{2}\left(
-c_{12}+z_{1}c_{22}\right) \right.  \nonumber \\
&&\left. +\left( z_{2}\right) ^{2}\left( -c_{11}+z_{1}c_{21}\right) \right]
/\left( z_{1}-z_{2}\right) ~,
\end{eqnarray}%
\begin{equation}
\alpha _{21}=\left[ -c_{15}+z_{1}c_{25}+z_{2}\left(
-c_{14}+z_{1}c_{24}\right) \right] /\left( z_{1}-z_{2}\right) ~,
\end{equation}%
\begin{equation}
\alpha _{20}=\left( -c_{16}+z_{1}c_{26}\right) /\left( z_{1}-z_{2}\right) ~.
\end{equation}

\textbf{Remark 4-1}. The $4$ \textit{constraints} (\ref{4ConstraintsFinal})
coincide with the formulas (\ref{3ConstraintsFinal}) and (\ref%
{3ThirdFourthCon}); the formulas (\ref{4xnt}) come from the eqs. (\ref{xnyn}%
), (\ref{3Alanda}), and (\ref{Cynt}); the formulas (\ref{4zn}) coincide with
the eqs. (\ref{3zn}); the formulas (\ref{4wnt}) come from the eqs. (\ref%
{Cynt}), (\ref{Cynpm}) and (\ref{2ynt}); the formulas (\ref{4w1w20}) come
from the eqs. (\ref{Cynt}), (\ref{y12x12}) and (\ref{3Alanda}); the formulas
(\ref{4wnpm}) come from the eqs. (\ref{Cynt}), (\ref{AalphaA}), (\ref%
{ynplusminusbetan}) and (\ref{3Alanda}); the $6$ formulas (\ref{4alphanl})
come from the eqs. (\ref{3a1L}), (\ref{3a2L}), (\ref{Calpha}), and (\ref%
{3Alanda}). $\blacksquare $

\textbf{Remark 4-2}. The special subcases of this system which feature the
remarkable property to be \textit{isochronous }are clearly those
characterized by the requirement that the $2$ parameters $\beta _{n}$---see (%
\ref{4wnpm}) with (\ref{4alphanl}) and (\ref{4zn})---be both \textit{rational%
} multiples of the \textit{same imaginary} number:
\end{subequations}
\begin{equation}
\beta _{n}=\mathbf{i}\rho _{n}\omega ~,~~~n=1,2~,
\end{equation}%
where $\mathbf{i}\omega $ is an \textit{arbitrary imaginary} number and $%
\rho _{n}$ are $2$ \textit{real} (positive or negative) \textit{rational}
numbers (this is rather obvious, but in case of doubt see, for instance,
\cite{C2008}). While, if one of the $2$ parameters $\beta _{n}$ is an
\textit{arbitrary purely imaginary} number and the other is \textit{not} a
\textit{purely imaginary number}, then the system (\ref{41}) is
\emph{asymptotically isochronous} (see \cite{CGU2008}). $\blacksquare $

\section{Two special cases of the system (\protect\ref{1})}

The treatment reported up to this point has assumed that the $12$
coefficients $c_{nj}$ in (\ref{1}) take \textit{generic} values (except, of
course, for satisfying the $4$ \textit{constraints} (\ref{4ConstraintsFinal}%
)). However in several context this is \textit{not} the case; for instance
the subclass of systems (\ref{1}) characterized by the restrictions $%
c_{13}=c_{15}=c_{21}=c_{24}=0$ is relevant in many \textit{applicable}
contexts; and the system (\ref{1}) with \textit{homogeneous} second-degree
polynomial right-hand sides---i. e., with $c_{nj}=0$ for $n=1,2$ and $%
j=4,5,6 $---also deserves a special treatment, in order to compare the
findings presented in the present paper with those reported in the recent
paper \cite{CCL2020}. These $2$ special cases are treated in the following $%
2 $ subsections of this \textbf{Section 5}.

\subsection{The subcase of (\protect\ref{1}) with $%
c_{13}=c_{15}=c_{21}=c_{24}=0$}

In many \textit{applicative} contexts it is unreasonable to assume that the
time-evolution of the dependent variable $x_{n}\left( t\right) $ is
influenced by agents (represented by terms in the right-hand sides of the
ODEs (\ref{1})) which depend \textit{only} on the other variable $%
x_{n+1}\left( t\right) $ (with $n=1,2~\mod[2]$). Hence the subcase of
the system (\ref{1}) characterized by the restrictions
\begin{subequations}
\label{cfgh}
\begin{equation}
c_{13}=c_{15}=c_{21}=c_{24}=0~  \label{vanishingCs}
\end{equation}%
deserves special attention, featuring indeed in many \textit{applicative}
contexts. For this reason in this \textbf{Subsection 5.1} we focus on the
special case of (\ref{1}) characterized by these limitations (\ref%
{vanishingCs}), introducing moreover---for notational simplicity---the
following new notation for the $8$ remaining coefficients $c_{nj}$ in (\ref%
{1}):
\begin{eqnarray}
c_{11} &=&f_{11}~,~~~c_{12}=f_{12}~,~~~c_{14}=g_{1}~,~~~c_{16}=h_{1}~,
\nonumber \\
c_{22} &=&f_{21}~,~~~c_{23}=f_{22}~,~~~c_{25}=g_{2}~,~~~c_{26}=h_{2}~;
\label{4cf}
\end{eqnarray}%
so that the system (\ref{1}) reads hereafter (in this \textbf{Subsection 5.1}%
) as follows:
\end{subequations}
\begin{subequations}
\label{4Syst1}
\begin{equation}
\dot{x}_{n}=x_{n}\left( f_{n1}x_{1}+f_{n2}x_{2}+g_{n}\right) +h_{n}~,~~~n=1,2
\label{4xndot}
\end{equation}%
namely%
\begin{equation}
\dot{x}_{1}=x_{1}\left( f_{11}x_{1}+f_{12}x_{2}+g_{1}\right) +h_{1}~,
\label{4x1dot}
\end{equation}%
\begin{equation}
\dot{x}_{2}=x_{2}\left( f_{21}x_{1}+f_{22}x_{2}+g_{2}\right) +h_{2}~.
\label{4x2dot}
\end{equation}

\textbf{Remark 5.1-1}. Clearly this system of $2$ coupled nonlinear ODEs is
invariant under the following transformation:
\end{subequations}
\begin{equation}
x_{1}\left( t\right) \Leftrightarrow x_{2}\left( t\right)
~,~~f_{11}\Leftrightarrow f_{22}~,~~f_{12}\Leftrightarrow
f_{21}~,~~g_{1}\Leftrightarrow g_{2}~,~~h_{1}\Leftrightarrow h_{2}~;
\label{4Transx1x2}
\end{equation}%
which clearly replaces the analogous transformation reported in \textbf{%
Remark 1-2}. $\blacksquare $

The interested reader will easily verify that a direct adaptation to the
present case of the final findings reported above (see \textbf{Section 4})
is \textit{a priori unjustified}, because the conditions (\ref{cfgh}) render
illegitimate some of the steps performed in that section---where the
treatment was indeed based on the assumption that the coefficients $c_{nj}$
in (\ref{1}) have \textit{generic} values (except for satisfying the
constraints (\ref{4ConstraintsFinal})). So below (in this \textbf{Subsection
5.1}) we review the above treatment, adapting it to the new situation. On
the other hand, for the convenience of the reader who is only interested in
the solution of the dynamical system (\ref{4Syst1}) and not in the details
of how that solution has been obtained, we report in the following \textbf{%
Subsubsection 5.1.1}---in analogy to what we did in \textbf{Section 4} for
the dynamical system (\ref{1})---the explicit solution of the system (\ref%
{4Syst1}); even at the cost of some repetitions.

The $12$ equations (\ref{c1j}) and (\ref{c2j}) read now (via (\ref%
{vanishingCs}) and (\ref{4cf})) as follows:
\begin{subequations}
\label{4fgh}
\begin{equation}
f_{11}=\left[ a_{12}A_{11}\left( A_{22}\right) ^{2}+a_{22}A_{12}\left(
A_{21}\right) ^{2}\right] /D^{2}~,  \label{4fgha}
\end{equation}%
\begin{equation}
f_{12}=-2A_{11}A_{12}\left[ a_{12}A_{22}+a_{22}A_{21}\right] /D^{2}~,
\label{4fghb}
\end{equation}%
\begin{equation}
A_{11}A_{12}\left( a_{12}A_{12}+a_{22}A_{11}\right) =0~,  \label{4fgh3}
\end{equation}%
\begin{equation}
g_{1}=\left( a_{11}A_{11}A_{22}-a_{21}A_{12}A_{21}\right) /D~,  \label{4fghd}
\end{equation}%
\begin{equation}
a_{11}=a_{21}~,  \label{4fghe}
\end{equation}%
\begin{equation}
h_{1}=a_{10}A_{11}+a_{20}A_{12}~;  \label{4fghf}
\end{equation}%
\begin{equation}
f_{22}=\left[ a_{22}A_{22}\left( A_{11}\right) ^{2}+a_{12}A_{21}\left(
A_{12}\right) ^{2}\right] /D^{2}~,  \label{4fghg}
\end{equation}%
\begin{equation}
f_{21}=-2A_{22}A_{21}\left( a_{22}A_{11}+a_{12}A_{12}\right) /D^{2}~,
\label{h4fgh}
\end{equation}%
\begin{equation}
A_{22}A_{21}\left( a_{22}A_{21}+a_{12}A_{22}\right) =0~,  \label{4fghi}
\end{equation}%
\begin{equation}
g_{2}=\left( a_{21}A_{11}A_{22}-a_{11}A_{12}A_{21}\right) /D~,  \label{4fghj}
\end{equation}%
\begin{equation}
h_{2}=a_{10}A_{21}+a_{20}A_{22}~,  \label{4fghk}
\end{equation}%
of course with $D$ defined as above, see (\ref{AA}).

\textbf{Remark 5.1-2}. In this \textbf{Subsection 5.1}, as mentioned above,
we assume that the parameter $D$ does \textit{not} vanish, but we do not
exclude the possibility that one of the $4$ parameters $A_{nm}$ vanish (in
contrast with \textbf{Remark 2-2}). $\blacksquare $

It is now again convenient to introduce the $2$ auxiliary parameters $%
z_{1}=A_{11}/A_{21}$ and $z_{2}=A_{12}/A_{22}$ (see (\ref{3z1z2})). Then the
$2$ eqs. (\ref{3Anm2}) and (\ref{3Anm1}) yield (via (\ref{cfgh}))
\end{subequations}
\begin{subequations}
\begin{equation}
\left[ \left( f_{21}-f_{11}\right) z_{n}+f_{22}-f_{12}\right]
z_{n}=0~,~~~n=1,2~,
\end{equation}%
implying
\end{subequations}
\begin{subequations}
\label{4z1z2}
\begin{equation}
z_{1}=0\,,~~~z_{2}=\frac{f_{22}-f_{12}}{f_{11}-f_{21}}~,  \label{4z1z2a}
\end{equation}%
or%
\begin{equation}
z_{2}=0\,,~~~z_{1}=\frac{f_{22}-f_{12}}{f_{11}-~f_{21}}~,  \label{4z1z2b}
\end{equation}%
since (see \textbf{Remark 5.1-2}) we exclude the solution $z_{1}=z_{2}$
which implies $D=0$ (see (\ref{AA}) and (\ref{3z1z2})). Note that---via (\ref%
{3z1z2})---$z_{1}=0$ implies $A_{11}=0$ and likewise $z_{2}=0$ implies $%
A_{12}=0$, each of these $2$ equalities reducing eq. (\ref{4fgh3}) to the
\textit{identity} $0=0$.

Let us now see how the remaining $11$ eqs. (\ref{4fgh}) simplify in the $%
z_{1}=0$ case, when clearly
\end{subequations}
\begin{subequations}
\label{4aAfgh}
\begin{equation}
A_{11}=0  \label{4A11D}
\end{equation}%
implying $D=-A_{12}A_{21}$; for analogous results in the $z_{2}=0$ case see
below \textbf{Remark 5.1.1-1}.

As already mentioned above, this eq. (\ref{4A11D}) implies that eq. (\ref%
{4fgh3}) holds identically ($0=0$); hence only the following $10$ equations
remain (note that they are reported below in a somewhat different order than
in (\ref{4fgh})):%
\begin{equation}
a_{22}A_{21}+a_{12}A_{22}=0~,  \label{4aaAA}
\end{equation}%
\begin{equation}
a_{11}=a_{21}~;
\end{equation}%
\begin{equation}
f_{11}=a_{22}/A_{12}~,  \label{4f11}
\end{equation}%
\begin{equation}
f_{12}=0~,  \label{4f12EqZero}
\end{equation}%
\begin{equation}
f_{21}=-2a_{12}A_{22}/\left( A_{12}A_{21}\right) ~,  \label{4f21}
\end{equation}%
\begin{equation}
f_{22}=a_{12}/A_{21}~;
\end{equation}%
\begin{equation}
g_{1}=a_{21}~,
\end{equation}%
\begin{equation}
g_{2}=a_{11}~;
\end{equation}%
\begin{equation}
h_{1}=a_{20}A_{12}~,
\end{equation}%
\begin{equation}
h_{2}=a_{10}A_{21}+a_{20}A_{22}~.
\end{equation}%
Clearly the first $3$ of these $11$ equations (\ref{4aAfgh}) provide $3$
\textit{constraints} on the $8$ parameters $A_{nm}$ and $a_{nm}$; while the
last $8$ of these $12$ eqs. (\ref{4aAfgh}) express \textit{explicitly---}in
terms of the $8$ parameters $A_{nm}$ and $a_{nm}$---the $8$ parameters $%
f_{nm},$ $g_{n},$ $h_{n}$ which characterize the system of ODEs (\ref{4Syst1}%
).

The next task is to invert the last $8$ formulas (\ref{4aAfgh}), namely to
express in terms of the $8$ parameters $f_{nm},$ $g_{n},$ $h_{n}$ featured
by the system (\ref{4Syst1}), the $8$ parameters $A_{nm},$ and $a_{nm}$
which characterize the explicit solution of this system (\ref{4Syst1}) via
the formulas of \textbf{Section 1} complemented by the restrictions and
redefinitions (\ref{cfgh}); and as well to identify---most importantly---the
\textit{constraints} implied by our treatment on the $8$ parameters $f_{nm},$
$g_{n},$ $h_{n}$.

These findings can be obtained by appropriately specializing the formulas
obtained in the preceding \textbf{Section 3}, i. e. by inserting in them the
formulas (\ref{cfgh}) as well as the findings reported above in this \textbf{%
Subsection 5.1}.

In this manner from the $3$ eqs. (\ref{3a1L}) we get (using the \textit{%
constraint} $f_{12}=0$ already obtained above, see (\ref{4f12EqZero}))
\end{subequations}
\begin{subequations}
\label{4a1L}
\begin{equation}
a_{12}=A_{21}f_{22}~,
\end{equation}%
\begin{equation}
a_{11}=g_{2}~,
\end{equation}%
\begin{equation}
a_{10}=\left( A_{12}h_{2}-A_{22}h_{1}\right) /\left( A_{12}A_{21}\right) ~;
\end{equation}%
and likewise from the $3$ eqs. (\ref{3a2L}) we get
\end{subequations}
\begin{subequations}
\label{4a2L}
\begin{equation}
a_{22}=A_{12}f_{11}~,
\end{equation}%
\begin{equation}
a_{21}=g_{1}~,
\end{equation}%
\begin{equation}
a_{20}=h_{1}/A_{12}~.
\end{equation}

Next, let us look at the $6$ eqs. (\ref{3Anmc1L}) and (\ref{Anm4}), using
again the \textit{constraint} $f_{12}=0$ (see (\ref{4f12EqZero})) to
simplify some of them.

The $3$ eqs. (\ref{3Anmc1L}) read then as follows:
\end{subequations}
\begin{subequations}
\begin{equation}
A_{12}f_{21}+2A_{22}f_{22}=0~,
\end{equation}%
\begin{equation}
A_{12}\left( f_{11}-f_{21}\right) -A_{22}f_{22}=0~,
\end{equation}%
\begin{equation}
g_{1}-g_{2}=0~.  \label{g1g2}
\end{equation}

It is now easily seen that the first $2$ of these $3$ equations imply (since
we assume that $A_{12}$ and $A_{22}$ do \textit{not} vanish) the following
\textit{second constraint} on the $2$ parameters $f_{11}$ and $f_{21}$:
\end{subequations}
\begin{equation}
f_{21}=2f_{11}~;  \label{4f21f11}
\end{equation}%
and we moreover find from eq. (\ref{g1g2}) the following \textit{third
constraint},%
\begin{equation}
g_{1}=g_{2}~.  \label{5g1g2}
\end{equation}

On the other hand the $3$ eqs. (\ref{Anm4}) are \textit{identically
satisfied---}i. e., $0=0$---thanks to (\ref{4A11D}), to the conditions (\ref%
{cfgh}) and, again, to the constraint $f_{12}=0$, see (\ref{4f12EqZero}).

So, let us summarize the findings in this $z_{1}=0$ case. There are the $2$
\textit{constraints} $f_{12}=0$ (see (\ref{4f12EqZero})) and $f_{21}=2f_{11}$
(see (\ref{4f21f11})) on the parameters $f_{nm}$ of the system (\ref{4Syst1}%
), and the third \textit{constraint} $g_{1}=g_{2}$; note that the first $2$
of these $3$ \textit{constraints} imply $z_{2}=-f_{22}/f_{11}$. Provided
these $3$ \textit{constraints} are satisfied, the \textit{explicit} solution
of the \textit{initial-values }problem for this system (\ref{4Syst1}) is
provided by the treatment of \textbf{Section 2}, of course\ with the
parameters $c_{nj}$ replaced by the parameters $f_{nm}$, $g_{n}$, $h_{n}$ as
implied by the relations (\ref{cfgh}), and with the parameters $a_{n\ell }$
expressed by the formulas (\ref{4a1L}) and (\ref{4a2L}) in terms of the $8$
parameters $f_{nm}$, $g_{n}$, $h_{n}$, and also of the $4$ parameters $%
A_{nm} $. As for these latter parameters, they are themselves determined in
terms of the parameters $f_{nm}$ as follows (the last of these formulas is
of course implied by $A_{12}=z_{2}A_{22}$ with $z_{2}=-f_{22}/f_{11}$, see
above):
\begin{equation}
A_{21}=\lambda _{1}~,~~~A_{22}=\lambda _{2}~,~~~A_{11}=0~,~~~A_{12}=-\left(
f_{22}/f_{11}\right) \lambda _{2}~.
\end{equation}%
Here $\lambda _{1}$ and $\lambda _{2}$ are again $2$ \textit{arbitrary}
(nonvanishing) parameters, which can be \textit{freely} assigned because
their values do not influence the solution of the problem (as explained at
the end of \textbf{Section 3} and in \textbf{Appendix C} in the context of
the more general case of the system (\ref{1}); and see also \textbf{%
Subsubsection 5.1.1}).

The corresponding solution of the \textit{initial-values} problem for the
system (\ref{4Syst1}) is reported in the following \textbf{Subsubsection
5.1.1}.

Finally, let us conclude this \textbf{Subsection 5.1} by emphasizing that
the $3$ \textit{constraints }(\ref{4f12EqZero}), (\ref{4f21f11}) and (\ref%
{5g1g2}) entail a significant limitation on the generality of the system (%
\ref{4Syst1}) treated in this section; for instance, they exclude models of
the Lotka-Volterra type, which require (at least!) an arbitrary (\textit{%
nonvanishing}) assignment of the parameter $f_{12}$. Nevertheless the fact
that the \textit{initial-values} problem for the system (\ref{4Syst1}) can
be explicitly solved provided only the $3$ restrictions indicated above hold
seems a significant new finding.

\subsubsection{Solution of the initial-values problem for the system (%
\protect\ref{4Syst1})}

In this \textbf{Subsubsection 5.1.1} we report the solution of the \textit{%
initial-values problem} for the system characterized by the following $2$
equations of motion:
\begin{subequations}
\label{511System}
\begin{equation}
\dot{x}_{1}=x_{1}\left( f_{1}x_{1}+g\right) +h_{1}~,
\end{equation}%
\begin{equation}
\dot{x}_{2}=x_{2}\left( 2f_{1}x_{1}+f_{2}x_{2}+g\right) +h_{2}~,
\end{equation}%
which correspond to the system (\ref{4Syst1}) treated in this \textbf{%
Section 5.1 }by taking into account the following \textit{constraints }and
\textit{simplified} notation:
\end{subequations}
\begin{subequations}
\label{511ConSimp}
\begin{equation}
c_{12}=c_{13}=c_{15}=c_{21}=c_{24}=0~,  \label{511Con}
\end{equation}
\begin{eqnarray}
c_{11} &=&f_{11}=f_{1}~,~~c_{14}=c_{25}=g_{1}=g_{2}=g~,~~c_{16}=h_{1}~,
\nonumber \\
c_{22}
&=&f_{21}=2f_{1}~,~~c_{23}=f_{22}=f_{2}~,~~c_{23}=f_{22}=f_{2}~;~~c_{26}=h_{2}~,
\label{511Simp}
\end{eqnarray}%
implied by our treatment, see above.

The solution is then provided by the following formulas (obtained via (\ref%
{511ConSimp}) from the corresponding solution reported in \textbf{Section 4}%
):
\end{subequations}
\begin{equation}
x_{1}\left( t\right) =-\left( f_{2}/f_{1}\right) \xi _{2}\left( t\right)
~,~~~x_{2}\left( t\right) =\xi _{1}\left( t\right) +\xi _{2}\left( t\right)
~;
\end{equation}%
\begin{subequations}
\begin{eqnarray}
\xi _{n}\left( t\right) &=&\frac{\xi _{n}\left( 0\right) \left[ \xi
_{n+}-\xi _{n-}\exp \left( \gamma _{n}t\right) \right] -\xi _{n+}\xi _{n-}%
\left[ 1-\exp \left( \gamma _{n}t\right) \right] }{\xi _{n+}\exp \left(
\gamma _{n}t\right) -\xi _{n-}+\xi _{n}\left( 0\right) \left[ 1-\exp \left(
\gamma _{n}t\right) \right] }~,  \nonumber \\
n &=&1,2~;  \label{511ksin}
\end{eqnarray}%
\begin{equation}
\xi _{1}\left( 0\right) =x_{2}\left( 0\right) +\left( f_{1}/f_{2}\right)
x_{1}\left( 0\right) ~,~~~\xi _{2}\left( 0\right) =-\left(
f_{1}/f_{2}\right) x_{1}\left( 0\right) ~,
\end{equation}%
\begin{equation}
\xi _{n\pm }=\left( -\eta _{n1}\pm \gamma _{n}\right) /\left( 2\eta
_{n2}\right) ~,~~~\gamma _{n}=\sqrt{\left( \eta _{n1}\right) ^{2}-4\eta
_{n0}\eta _{n2}}~,~~n=1,2~;
\end{equation}%
\end{subequations}
\begin{subequations}
\begin{equation}
\eta _{12}=f_{2}~,~~~\eta _{11}=g~,~~~\eta _{10}=\left( f_{1}/f_{2}\right)
h_{1}+h_{2}~,
\end{equation}

\begin{equation}
\eta _{22}=-f_{2}~,~~~\eta _{21}=g~,~~~\eta _{20}=-\left( f_{1}/f_{2}\right)
h_{1}~;
\end{equation}%
implying
\end{subequations}
\begin{equation}
\gamma _{1}=\sqrt{g^{2}-4f_{2}\eta _{10}}~,~~~\gamma _{2}=\sqrt{%
g^{2}+4f_{2}\eta _{20}}~.
\end{equation}

\textbf{Remark 5.1.1-1}. The diligent reader will check that the \textit{%
same }result is obtained in the alternative case with $z_{2}=0$ (i. e. with (%
\ref{4z1z2b}) replacing (\ref{4z1z2a})); of course provided all the
corresponding notational changes are made, consisting essentially to an
exchange of the roles of the auxiliary variables $\xi _{1}\left( t\right) $
and $\xi _{2}\left( t\right) $---themselves corresponding, of course up to
an appropriate rescaling, to the variables $y_{1}\left( t\right) $ and $%
y_{2}\left( t\right) $, in analogy to the treatment detailed, for the more
general system (\ref{1}), in \textbf{Appendix C}. $\blacksquare $

\subsection{The subcase of (\protect\ref{1}) with homogeneous second-degree
polynomial right-hand sides, i. e. $c_{nj}=$ $0$ for $n=1,2$ and $j=4,5,6$}

The special case of the dynamical system (\ref{1}) treated in this \textbf{%
Subsection 5.2} is characterized by the following ODEs:
\begin{subequations}
\label{3CCL}
\begin{equation}
\dot{x}_{n}\left( t\right) =c_{n1}\left[ x_{1}\left( t\right) \right]
^{2}+c_{n2}x_{1}\left( t\right) x_{2}\left( t\right) +c_{n3}\left[
x_{2}\left( t\right) \right] ^{2}~,~~~n=1,2~,  \label{3CCLa}
\end{equation}%
namely%
\begin{equation}
\dot{x}_{1}\left( t\right) =c_{11}\left[ x_{1}\left( t\right) \right]
^{2}+c_{12}x_{1}\left( t\right) x_{2}\left( t\right) +c_{13}\left[
x_{2}\left( t\right) \right] ^{2}~,  \label{3CCLb}
\end{equation}%
\begin{equation}
\dot{x}_{2}\left( t\right) =c_{21}\left[ x_{1}\left( t\right) \right]
^{2}+c_{22}x_{1}\left( t\right) x_{2}\left( t\right) +c_{23}\left[
x_{2}\left( t\right) \right] ^{2}~.  \label{3CCLc}
\end{equation}%
A large class of subcases of this system---identified by explicit
restrictions on the $6$ parameters $c_{nk}$ ($n=1,2;$ $k=1,2,3$) and
characterized by the property to be \textit{solvable} by \textit{algebraic}
operations---has been identified in the recent paper \cite{CCL2020}. In this
\textbf{Subsection 5.2} we compare the results obtained in this paper \cite%
{CCL2020} \ with those obtained in the present paper. The main conclusion of
this comparison is that there is, of course, a certain overlap among the
cases treated in the present paper and those treated in \cite{CCL2020};
however subcases of (\ref{1}) identified as \textit{solvable} in \cite%
{CCL2020} are \textit{not} included in the treatment provided in the present
paper; and likewise subcases of (\ref{1}) identified as \textit{solvable} in
the present paper are \textit{not} included in the treatment provided by
\cite{CCL2020}. Hence the $2$ approaches ---with their similarities and
their differences---are in some sense \textit{complementary}. This is
explained in detail in the following $2$ Subsubsections.

\textbf{Remark 5.2-1}. There is a significant difference among the
treatments of \cite{CCL2020} and the present paper. In \cite{CCL2020} the
systems identified are \textit{algebraically solvable} in the following
sense: the computation of their $t$-evolution is reduced to the evaluation
of the zeros of a $t$-dependent polynomial $P_{N}\left( x;t\right) $ of
(finite) order $N$ in its argument $x$ ($N$ being a \textit{positive integer}
whose value depends on the particular model under consideration), the $t$%
-evolution of which is \textit{explicitly} known (while the expressions of
the solutions $x_{n}\left( t\right) $ of the system (\ref{3CCL}) can of
course only be obtained \textit{explicitly }if $N\leq 4$). In the present
paper the systems identified as \textit{solvable} allow the \textit{explicit}
display of the $t$-evolution $x_{n}\left( t\right) $ of the solutions of
their initial-values problem, as reported in \textbf{Section 4}.~$%
\blacksquare $

\textbf{Remark 5.2-2}. Another significant difference among the class of
systems treated in \cite{CCL2020} and in the present paper is that---while
the system (\ref{1}) features the $12$ \textit{a priori arbitrary }%
coefficients $c_{nj}$---the solutions obtained in \cite{CCL2020} generally
feature $9$ ($9=2+6+1$: see below \textbf{Subsubsection 5.2.2}) \textit{%
freely assigned} parameters, while those obtained in the present paper
feature $8$ ($8=12-4$: see above \textbf{Section 4}) \textit{freely assigned}
parameters (in addition of course, in both cases, to the $2$ initial data $%
x_{n}\left( 0\right) $). But, as shown in the following $2$ \textbf{%
Subsubsections }(see their titles), neither one of these two subclasses of
the system (\ref{1})---that treated in \cite{CCL2020} and that treated in
the present paper---includes the other subclass: a confirmation that these $%
2 $ papers are actually \textit{complementary}. $\blacksquare $

\subsubsection{Subcases of (\protect\ref{1}) shown to be \textit{solvable}
in \protect\cite{CCL2020} which are \textit{not} included among those shown
to be \textit{solvable} in the present paper}

In \cite{CCL2020} it is noted that essentially the entire class of systems (%
\ref{3CCL}) can be reduced to the following simpler system (see eq. (6) of
\cite{CCL2020}):
\end{subequations}
\begin{equation}
\dot{x}_{1}=x_{1}x_{2}~,~~~\dot{x}_{2}=A\left[ \left( x_{1}\right)
^{2}+\left( x_{2}\right) ^{2}\right] +Bx_{1}x_{2}~,  \label{31CCL}
\end{equation}%
and that this system can be solved by \textit{algebraic} operations if the $%
2 $ parameters $A$ and $B$ are suitably restricted, for instance \textit{%
sufficient} conditions are that
\begin{equation}
A=\frac{n+q-1}{n+q}~,~~~B=\pm \frac{n-q}{n+q}\sqrt{\frac{n+q-1}{nq}}
\label{3ABnq}
\end{equation}%
with $n$ an \textit{arbitrary positive integer}, and $q$ an \textit{arbitrary%
}, possibly \textit{complex}, \textit{rational} number (see eqs. (17) of
\cite{CCL2020}; there also are other possibilities, see eqs. (18-20) of \cite%
{CCL2020}, but we do not need to evoke them to make our point).

On the other hand it is easily seen that the system (\ref{31CCL}), which
corresponds to the system (\ref{1}) only if \textit{all} the $12$ parameters
$c_{nj}$ vanish except for the following $4$ of them,%
\begin{equation}
c_{12}=1~,~~~c_{21}=A~,~~~c_{22}=B~,~~~c_{23}=A~,  \label{3cAB}
\end{equation}%
entails, via the $4$ constraints (\ref{4ConstraintsFinal}), either%
\begin{equation}
A=B=0~,
\end{equation}%
which is also consistent with (\ref{3ABnq}) (say, with $q=1-n)$, or%
\begin{equation}
B=0~,~~~A=1/2~,
\end{equation}%
which is consistent (say, with $n=q=1$); both assignments, of course, much
less general than (\ref{3ABnq}) with $n$ an \textit{arbitrary positive
integer}, and $q$ an \textit{arbitrary}, possibly \textit{complex}, \textit{%
rational} number.

This shows that there are some special subcases of (\ref{1}) which are
\textit{solvable} both via the technique of \cite{CCL2020} and via the
technique of the present paper; and many more which are \textit{solvable}
via the technique of \cite{CCL2020} but are \textit{not solvable} via the
technique of the present paper. Q. E. D.

\subsubsection{Subcases of (\protect\ref{1}) shown to be \textit{solvable}
in the present paper which are \textit{not} included among those shown to be
\textit{solvable} in \protect\cite{CCL2020}}

Since the system (\ref{1})---even with the \textit{constraints} (\ref%
{4ConstraintsFinal})---is clearly more general than the system (\ref{3CCL})
treated in \cite{CCL2020}---because of the additional $6$ terms featuring
the coefficients $c_{nj}$ with $n=1,2$ and $j=4,5,6$, it might seem that
what we want to\ prove in this \textbf{Subsubsection 5.2.2} (see its title)
is altogether obvious. But the subclass of the system (\ref{1}) treated in
\cite{CCL2020} includes the \textit{additional} possibility to perform a
\textit{linear} transformation---with \textit{arbitrary time-independent}
coefficients---of the dependent variables. So this argument is \textit{not}
cogent.

But such a transformation---which generally features $6$ \textit{free }%
parameters---cannot change the dependence on the independent variable $t$
from \textit{algebraic} to \textit{exponential}; while the dependence on the
variable $t$ of the solutions reported in the present paper is indeed
generally \textit{exponential} (see for instance above eq. (\ref{4wnt})).

However this argument is still not entirely conclusive because of the
possibility to extend the system (\ref{3CCL}) via the simple \textit{%
invertible} change of dependent and independent variables%
\begin{equation}
x_{n}\left( t\right) =\exp \left( \lambda t\right) \zeta _{n}\left( \tau
\right) ~,~~~\tau =\left[ \exp \left( \lambda t\right) -1\right] /\lambda ~,
\label{3transttau}
\end{equation}%
which transforms the following system for $\zeta _{n}\left( \tau \right) ,$
reading
\begin{equation}
d\zeta _{n}\left( \tau \right) /d\tau =c_{n1}\left[ \zeta _{1}\left( \tau
\right) \right] ^{2}+c_{n2}\zeta _{1}\left( \tau \right) \xi _{2}\left( \tau
\right) +c_{n3}\left[ \zeta _{2}\left( \tau \right) \right] ^{2}~,~~~n=1,2~,
\end{equation}%
hence being included among those treated in \cite{CCL2020}, into the, also
\textit{autonomous}---and as well \textit{solvable }(via (\ref{3transttau}%
))---system%
\begin{equation}
\dot{x}_{n}\left( t\right) =\lambda x_{n}\left( t\right) +c_{n1}\left[
x_{1}\left( t\right) \right] ^{2}+c_{n2}x_{1}\left( t\right) x_{2}\left(
t\right) +c_{n3}\left[ x_{2}\left( t\right) \right] ^{2}~,~~~n=1,2~.
\end{equation}%
This new system features the \textit{additional free} parameter $\lambda $
and---most importantly with respect to the previous argument---its solutions
clearly feature now an \textit{exponential} dependence on the independent
variable $t$ (see (\ref{3transttau})).

But this implies that the solutions of this model---even after a linear
reshuffling of the dependent variables---can only feature a dependence on
the single exponential $\exp \left( \lambda t\right) $; while the solutions
obtained in the present paper feature the $2$, generally \textit{different},
exponentials $\exp \left( \beta _{n}t\right) $, $n=1,2$, see above eqs. (\ref%
{4wnt}) and (\ref{4wnpm}).

It is thereby shown that there indeed are subcases of (\ref{1}) shown to be
\textit{solvable} in the present paper which are \textit{not} included among
those shown to be \textit{solvable} in \cite{CCL2020}. Q. E. D.

\section{Conclusions and outlook}

The prototypical system (\ref{1}) and its subcases treated above in \textbf{%
Section 5 }have been investigated over time by top mathematicians and
subtend an enormous number of applied mathematics models in several
scientific fields. It is our hope to obtain analogous results for analogous
models in the future; for one such result see \cite{CP2020}.

\section{Acknowledgements}

The results reported in this paper have been mainly obtained by a
collaboration at a distance among its two authors (essentially via e-mails).
We would like to acknowledge with thanks $2$ grants, which shall facilitate
our future collaboration by allowing FP to visit (hopefully more than once)
in 2021 the Department of Physics of the University of Rome "La Sapienza":
one granted by that University, and one granted jointly by the Istituto
Nazionale di Alta Matematica (INdAM) of that University and by the
International Institute of Theoretical Physics (ICTP) in Trieste in the
framework of the ICTP-INdAM "Research in Pairs" Programme. Finally, we
gratefully acknowledge a special contribution by Fran\c{c}ois Leyvraz, who
pointed out a serious flaw in a preliminary version of this paper, the
elimination of which also entailed a substantial simplification of its
presentation.

\section{Appendix A}

In this \textbf{Appendix A} we tersely demonstrate the following elementary
fact, which clearly implies the result (\ref{ynt}): that the solution of the
\textit{initial-values} problem for the ODE
\begin{subequations}
\begin{equation}
\dot{y}\left( t\right) =a_{2}\left[ y\left( t\right) \right]
^{2}+a_{1}y\left( t\right) +a_{0}~,  \label{Aydot}
\end{equation}%
is provided by the formula%
\begin{equation}
y\left( t\right) =\frac{y_{+}\left[ y\left( 0\right) -y_{-}\right] -y_{-}%
\left[ y\left( 0\right) -y_{+}\right] \exp \left( \beta t\right) }{y\left(
0\right) -y_{-}-\left[ y\left( 0\right) -y_{+}\right] \exp \left( \beta
t\right) }~,  \label{Ayt}
\end{equation}%
with $y_{\pm }$ defined as follows:%
\begin{equation}
y_{\pm }=\left( -a_{1}\pm \beta \right) /\left( 2a_{2}\right) ~,~~~\beta =%
\sqrt{\left( a_{1}\right) ^{2}-4a_{0}a_{2}}~.  \label{Aypluminusbeta}
\end{equation}

Indeed the ODE (\ref{Aydot}) can clearly be reformulated as follows:
\end{subequations}
\begin{subequations}
\begin{equation}
\dot{y}=a_{2}\left( y-y_{+}\right) \left( y-y_{-}\right)
\end{equation}%
with $y_{\pm }$ defined by (\ref{Aypluminusbeta}); and then (again, via (\ref%
{Aypluminusbeta})) this ODE can be rewritten as follows:
\begin{equation}
\dot{y}\left[ \left( y-y_{+}\right) ^{-1}-\left( y-y_{-}\right) ^{-1}\right]
=\beta ~.
\end{equation}%
The integration of this ODE for the dependent variable $y\left( t^{\prime
}\right) $ over the independent variable $t^{\prime }$---from $t^{\prime }=0$
to $t^{\prime }=t$---clearly yields
\end{subequations}
\begin{equation}
\ln \left[ \frac{y\left( t\right) -y_{+}}{y\left( 0\right) -y_{+}}\right]
-\ln \left[ \frac{y\left( t\right) -y_{-}}{y\left( 0\right) -y_{-}}\right]
=\beta t~,
\end{equation}%
which coincides---after exponentiation---with (\ref{Ayt}). Q. E. D.

\section{Appendix B}

In this \textbf{Appendix B} we tersely outline the derivation of the
expressions (\ref{c1j}) and (\ref{c2j}) of the $12$ parameters $c_{nj}$ in
terms of the $10$ parameters $A_{nm}$ and $a_{n\ell }$.

The \textbf{first step} is to invert the relations (\ref{xnyn}), getting
\begin{subequations}
\label{y12x12}
\begin{equation}
y_{1}\left( t\right) =\left[ A_{22}x_{1}\left( t\right) -A_{12}x_{2}\left(
t\right) \right] /D~,  \label{y1xn}
\end{equation}%
\begin{equation}
y_{2}\left( t\right) =\left[ A_{11}x_{2}\left( t\right) -A_{21}x_{1}\left(
t\right) \right] /D~,  \label{y2xn}
\end{equation}%
where the quantity $D$ is defined as above, see (\ref{AA}).

The \textbf{second step} is to note that the relations (\ref{xnyn}) imply
\end{subequations}
\begin{subequations}
\begin{equation}
\dot{x}_{n}=A_{n1}\dot{y}_{1}+A_{n2}\dot{y}_{2}~,~~~n=1,2~,
\end{equation}%
hence, via the ODEs (\ref{yndot}),
\begin{eqnarray}
\dot{x}_{n}=A_{n1}\left[ a_{12}\left( y_{1}\right) ^{2}+a_{11}y_{1}+a_{10}%
\right] &&  \nonumber \\
+A_{n2}\left[ a_{22}\left( y_{2}\right) ^{2}+a_{21}y_{2}+a_{20}\right]
~,~~~n=1,2~. &&
\end{eqnarray}

The \textbf{third and last step} is to insert the expressions (\ref{y12x12})
of $y_{1}$ and $y_{2}$ in terms of $x_{1}$ and $x_{2}$ in the right-hand
sides of these ODEs. Then, via a bit of trivial if tedious algebra, there
obtains the system (\ref{1}) with the definitions (\ref{c1j}) and (\ref{c2j}%
) of the $12$ coefficients $c_{nj}$. Q. E. D.

\section{Appendix C}

In this \textbf{Appendix C} we show that the solutions $x_{n}\left( t\right)
$ of the dynamical system (\ref{1})---as treated above, see \textbf{Sections
2} and \textbf{3}---do \textit{not} depend on the \textit{free} parameters $%
\lambda _{n}$ introduced via the positions (\ref{3Alanda}).

To this end we insert the expressions (\ref{3Alanda}) of the parameters $%
A_{nm}$ in terms of the free parameters $\lambda _{n}$ in the formulas (\ref%
{3a1L}) and (\ref{3a2L}) expressing the $6$ parameters $a_{nk}$; in order to
display the very simple dependence of these $8$ parameters from the $2$ free
parameters $\lambda _{n}$. We thus easily find the following formulas:
\end{subequations}
\begin{subequations}
\label{AalphaA}
\begin{equation}
a_{n\ell }\equiv \left( \lambda _{n}\right) ^{\ell -1}\alpha _{n\ell }\left(
C\right) ~,~~~n=1,2~,~~~\ell =0,1,2~,  \label{Calpha}
\end{equation}%
where the notation $C$ indicates---above and hereafter---the set of the $12$
parameters $c_{nj}$, and the $6$ functions $\alpha _{n\ell }\left( C\right) $
are \textit{explicitly} displayed in \textbf{Section 4}, see (\ref{4alphanl}%
); of course to this end we also used the definitions (\ref{3zn}) of the $2$
auxiliary parameters $z_{n}$ in terms of the $4$ coefficients $c_{nm}$ ($%
n=1,2;$ $m=1,2$).

The next step is to insert the formulas (\ref{Calpha}) in the expressions (%
\ref{ynt}), getting thereby%
\begin{equation}
y_{n\pm }\equiv \left( \lambda _{n}\right) ^{-1}w_{n\pm }\left( C\right)
~,~~~\beta _{n}\equiv \beta _{n}\left( C\right) ~,~~~n=1,2~,  \label{Cynpm}
\end{equation}%
again with the functions $w_{n\pm }\left( C\right) $ and $\beta _{n}\left(
C\right) $ \textit{explicitly} displayed in \textbf{Section 4}, see (\ref%
{4wnt}) and (\ref{4wnpm}).

The insertion of these formulas in the expressions (\ref{2ynt}) of the
solutions $y_{n}\left( t\right) $ of the auxiliary dynamical system (\ref%
{yndot}) evidences the following, very simple, dependence of these functions
from the \textit{free} parameters $\lambda _{n}$:
\end{subequations}
\begin{equation}
y_{n}\left( t\right) \equiv \left( \lambda _{n}\right) ^{-1}w_{n}\left(
C,t\right) ~,~~~n=1,2~,  \label{Cynt}
\end{equation}%
where again the $2$ functions $w_{n}\left( C,t\right) $ are \textit{%
explicitly} displayed in \textbf{Section 4}, see (\ref{4wnt}).

And via the insertion in the expressions (\ref{xnyn}) of $x_{n}\left(
t\right) $ of these formulas (\ref{Cynt}), together with the expressions (%
\ref{3Alanda}) of $A_{nm}$, we conclude that the solutions $x_{n}\left(
t\right) $ are independent of the \textit{free} parameters $\lambda _{n}$;
as indeed displayed in \textbf{Section 4}, see the set of eqs. from (\ref%
{4xnt}) to (\ref{4alphanl}). Q. E. D.

\end{document}